\UseRawInputEncoding
\documentclass[12pt, reqno]{amsart}
\usepackage{amssymb,latexsym,amsmath,amscd,amsfonts}
\usepackage{latexsym}
\usepackage{mathtools}
\usepackage[mathscr]{eucal}
\usepackage{bm}

\newcommand{\ft}{\mathcal F_O}

\def ~{\hspace{1mm}}

\def\textmatrix#1&#2\\#3&#4\\{\bigl({#1 \atop #3}\ {#2 \atop #4}\bigr)}
\def\dispmatrix#1&#2\\#3&#4\\{\left({#1 \atop #3}\ {#2 \atop #4}\right)}
\newcommand{\beg}{\begin{equation}}
\newcommand{\eeg}{\end{equation}}
\newcommand{\ben}{\begin{eqnarray*}}
\newcommand{\een}{\end{eqnarray*}}

\newtheorem{thm}{Theorem}[section]

\newtheorem{lem}[thm]{Lemma}

\numberwithin{equation}{section} \theoremstyle{definition}
\newtheorem{defn}[thm]{Definition}

\def\textmatrix#1&#2\\#3&#4\\{\bigl({#1 \atop #3}\ {#2 \atop #4}\bigr)}
\def\dispmatrix#1&#2\\#3&#4\\{\left({#1 \atop #3}\ {#2 \atop #4}\right)}

\begin{document}
\title[Automorphism and the fundamental operator]
{Automorphisms and the fundamental operators associated with the
symmetrized tridisc}

\author[Bappa Bisai]{Bappa Bisai}
\address[Bappa Bisai]{Mathematics Department, Indian Institute of Technology Bombay,
Powai, Mumbai - 400076, India.} \email{bisai@math.iitb.ac.in}

\author[Sourav Pal]{Sourav Pal}
\address[Sourav Pal]{Mathematics Department, Indian Institute of Technology Bombay,
Powai, Mumbai - 400076, India.} \email{sourav@math.iitb.ac.in}



\thanks{The first by named author is supported by a Ph.D fellowship
of the University Grants Commissoin (UGC). The second named author
is supported by the Seed Grant of IIT Bombay, the CPDA and the
INSPIRE Faculty Award (Award No. DST/INSPIRE/04/2014/001462) of
DST, India.}

\begin{abstract}
The automorphisms of the symmetrized polydisc $\mathbb G_n$ are well-known and are given in the coordinates of the polydisc in \cite{E:Z}. We find an explicit formula for the automorphisms of $\mathbb G_n$ in its own coordinates. If
$\tau$ is an automorphism of $\mathbb G_n$, then
$\tau(S_1,\dots,S_{n-1},P)$ is a $\Gamma_n$-contraction, where a
$\Gamma_n$-contraction is a commuting $n$-tuple of Hilbert space
operators for which the closed symmetrized polydisc $\Gamma_n$ is
a spectral set. Corresponding to every $\Gamma_n$-contraction
$(S_1,\dots,S_{n-1},P)$, there exist $n-1$ unique operators
$A_1,\dots,A_{n-1}$ such that
\[
S_i-S_{n-i}^*P=D_PA_iD_P\,, \quad D_P=(I-P^*P)^{1/2}\,,
\]
for $i=1,\dots, n-1$. This unique $(n-1)$-tuple
$(A_1,\dots,A_{n-1})$, which is called the fundamental operator
tuple or $\ft$-tuple of $(S_1,\dots,S_{n-1},P)$ in literature,
plays central role in every section of operator theory on
$\Gamma_n$. We find an explicit form of the $\ft$-tuple of $\tau
(S_1,\dots,S_{n-1},P)$ when $n=3$. We show by an example that a $\Gamma_n$-contraction
may not have commuting $\ft$-tuple. Also, we obtain a necessary and sufficient condition under which
two $\Gamma_n$-contractions are unitarily equivalent.
\end{abstract}

\maketitle


\section{Introduction}

\noindent The open and closed \textit{symmetrized polydisc} or,
\textit{symmetrized $n$-disc} for $n\geq 2$, are the following
subsets of $\mathbb C^n$:
\begin{align*}
\mathbb G_n &=\left\{ \left(\sum_{1\leq i\leq n} z_i,\sum_{1\leq
i<j\leq n}z_iz_j,\dots,
\prod_{i=1}^n z_i \right): \,|z_i|< 1, i=1,\dots,n \right \}, \\
\Gamma_n & =\left\{ \left(\sum_{1\leq i\leq n} z_i,\sum_{1\leq
i<j\leq n}z_iz_j,\dots, \prod_{i=1}^n z_i \right): \,|z_i|\leq 1,
i=1,\dots,n \right \}.
\end{align*}
The symmetrized polydisc has attracted considerable interests in
recent decades because of its complex geometry \cite{costara,
E:Z}, rich function theory \cite{jarnicki, costara1, NiPfZw, A:Y,
A:L:Y} and associated beautiful operator theory
\cite{tirtha-sourav, tirtha-sourav1, S:B, sourav, spal12, spal1, Bisai-Pal1, S:P, B:P}. Numerous interesting articles have been written on the
symmetrized polydisc and only a few among them are mentioned here.
More could be found in the reference lists of the articles cited
here. Studying operators associated with a domain is always of
interest in its own right. In this article, we describe an
interrelation between automorphisms of the symmetrized polydisc
and the tuples of commuting operators that have the closed
symmetrized polydisc as a spectral set.

\begin{defn}
A tuple of $n$ commuting operators $(S_1,\dots,S_{n-1},P)$,
defined on a Hilbert space $\mathcal H$, for which $\Gamma_n$ is a
spectral set is called a $\Gamma_n$-contraction.
\end{defn}

\noindent An automorphism of a domain $\Omega \subseteq \mathbb
C^n$ is a bijective and bi-holomorphic self-map of $\Omega$. It is
well-known that any automorphism $\tau$ of $\mathbb{G}_n$ (see
\cite{E:Z}) is of the following form:
\begin{align}\label{eqn:01}
& \tau\left(\sum_{i=1}^{n}z_i,\sum_{1 \leq i<j\leq n}z_iz_j, \dots ,\prod_{i=1}^{n}z_i\right) \notag\\
 = &\left(\sum_{i=1}^{n}m(z_i),\sum_{1 \leq i<j\leq n}m(z_i)m(z_j), \dots , \prod_{i=1}^{n}m(z_i)\right),
\end{align}
where $z_1, \dots , z_n$ are in the unit disk $\mathbb{D}$ and $m$
is an automorphism of $\mathbb{D}$. An
automorphism $m$ of the unit disc $\mathbb D$ is a map
\begin{equation}\label{eqn:02}
m(z):= \beta\dfrac{z - a}{1 - \bar{a}z}\;, \text{ for some } a \in
\mathbb{D} \text{ and } \beta \in \mathbb T,
\end{equation}
where $\mathbb T$ is the unit circle in the complex plane. For
obvious reason, an automorphism of $\mathbb G_n$ is often denoted
by $\tau_m$. In Lemma \ref{lem1}, we find
an explicit form of an automorphism of $\mathbb G_n$ in its own coordinates.
As a consequence of this result, we obtain a necessary and sufficient condition under which two $\Gamma_n$-contractions become unitarily equivalent. This is given in Theorem \ref{thm:equiv}.\\

Also if $\tau$ is a $\mathbb{C}^n$-valued holomorphic map in a
neighbourhood $N(\Gamma_n)$ of $\Gamma_n$ that maps $\Gamma_n$
into itself, then by functional calculus $ \tau(S_1, \dots
,S_{n-1},P)$ makes sense for any $\Gamma_n$-contraction
$(S_1,\dots, S_{n-1},P)$. Indeed, if we denote
\[
(S_{1\tau}, \dots
,S_{(n-1)\tau},P_{\tau})= \tau(S_1, \dots ,S_{n-1},P),
\]
then $(S_{1\tau}, \dots ,S_{(n-1)\tau},P_{\tau})$ is also a
$\Gamma_n$-contraction. It was shown in \cite{S:P} that
corresponding to every $\Gamma_n$-contraction $(S_1,\dots,
S_{n-1},P)$ there exists a unique $(n-1)$-tuple $(A_1,\dots,
A_{n-1})$ such that
\[
S_i-S_{n-i}^*P=D_PA_iD_P\,,
\]
where $D_P=(I-P^*P)^{1/2}$. The $(n-1)$-tuple
$(A_1,\dots,A_{n-1})$ is called the \textit{fundamental operator
tuple} or in short form the $\ft$-\text{tuple} of the
$\Gamma_n$-contraction $(S_1,\dots,S_{n-1},P)$. The reason behind
carrying such a name is that it plays pivotal role in all aspects
of operator theory on $\Gamma_n\;,$ e.g., \cite{tirtha-sourav,
tirtha-sourav1, sourav, pal-shalit, spal1, S:P, J:S}. The main aim
of this paper is to find an explicit form of the $\ft$-tuple of
$(S_{1\tau}, \dots ,S_{(n-1)\tau},P_{\tau})$ in terms of the
$\ft$-tuple of $(S_1,\dots, S_{n-1},P)$. To avoid excessive
complexities in the expressions, we restrict our attention to
$n=3$ here. We believe that analogous expressions could be
obtained for an arbitrary $n$ in a similar fashion, only one needs
to deal with a bit more complicated calculations. The same program
for $n=2$ was carried out as a part of the paper
\cite{tirtha-sourav1}. Since substantial differences have been
witnessed between the operator theory on $\Gamma_2$ and
$\Gamma_3$, we find it worth to determine the $\ft$-pair of the
$\Gamma_3$-contraction $(S_{1\tau}, S_{2\tau},P_{\tau})$ which is
described in Theorem \ref{thm13}.\\

If $(s_1,\dots, s_{n-1},p)$ is a scalar $\Gamma_n$-contraction, that is, a point in $\Gamma_n$, then the $\ft$-tuple of $(s_1,\dots , s_{n-1},p)$ is a point in $\Gamma_{n-1}$. This was obtained by Costara in \cite{costara1} in a different format. So, a natural question arises when we deal with $\ft$-tuple of a $\Gamma_n$-contraction: is the $\ft$-tuple  of a $\Gamma_n$-contraction a $\Gamma_{n-1}$-contraction ? We shall show that the answer is negative. Indeed, in the last section we shall produce a $\Gamma_n$-contraction whose $\ft$-tuple is not even commutative.

\section{Automorphisms of $\mathbb G_n$ in its own coordinates}

\noindent The automorphisms of the symmetrized polydisc $\mathbb
G_n$ were explicitly determined in \cite{E:Z}. Here we determine
the (same) formula in the coordinates of $\mathbb G_n$.

\begin{lem}\label{lem1}
    Suppose $(s_1, \dots, s_{n-1},s_n) \in \Gamma_n$ and $\tau$ is
    an automorphism of $\mathbb{G}_n$. Then $\tau(s_1, \dots, s_{n-1},s_n)
    = (s_{1\tau}, \dots, s_{(n-1)\tau},s_{n\tau})$, where
    \begingroup
    \allowdisplaybreaks
    \begin{align*}
    & s_{i\tau}\\
    &=   \beta^i\; \dfrac{\sum\limits_{j=1}^{i-1}(-a)^{i-j}
    [\sum\limits_{k=0}^j\{{{n-j}\choose{i-j+k}}{{j}\choose{k}}|a|^{2k}\}s_j]
    + \sum\limits_{j = i}^n(-\bar{a})^{j-i}[\sum\limits_{k=0}^i{{j}\choose {i-k}}{{n-j}\choose {k}}|a|^{2k}]s_j}{1
    + \sum\limits_{i=1}^n(-1)^i(\bar{a})^is_i} \\
    & + \beta^i \;\dfrac{(-1)^i{{n}\choose{i}}a^i}{1
    + \sum\limits_{i=1}^n(-1)^i(\bar{a})^is_i}, \text{ for any } a \in \mathbb{D} \text{ and }\beta \in \mathbb{T}.
    \end{align*}
    \endgroup
    \begin{proof}
        We prove by induction on $n$. Clearly, the lemma is
        true for $ n=2 $. Suppose the lemma is true for $n$.
        Consider $(s_1, \dots, s_{n+1}) \in \Gamma_{n+1}$. Then there exists
        $(z_1, \dots, z_{n+1}) \in \overline{\mathbb{D}}^{n+1}$ such that
        $\pi_{n+1}(z_1, \dots, z_{n+1}) = (s_1, \dots, s_{n+1})$.
        Let $\pi_n(z_1, \dots, z_n) = (s_1^{\prime}, \dots, s_n^{\prime})$. Then clearly,
        \begingroup
        \allowdisplaybreaks
        \begin{align*}
            & s_1 = s_1^{\prime} + z_{n+1}\\
            &s_i = s_i^{\prime} + s_{i-1}^{\prime}z_{n+1}, \text{ for } 1<i<n+1\\
            &s_{n+1} = s_n^{\prime}z_{n+1}.
        \end{align*}
        \endgroup
        Suppose $\tau_{n+1}\in \text{Aut}(\mathbb{G}_{n+1})$
        and

\begin{align*}
        \tau_{n+1}(s_1, \dots, s_{n+1}) & =(s_{1\tau},\dots,
        s_{(n+1)\tau})\\&
        = \pi_{n+1}\left(\beta\dfrac{z_1 - a}{1 - \bar{a}z_1},
        \dots, \beta\dfrac{z_{n+1} - a}{1 -
        \bar{a}z_{n+1}}\right)\,,
\end{align*}
        for some $a \in \mathbb{D}$ and $\beta \in \mathbb{T}$.
        Let $\tau_n \in \text{Aut}(\mathbb{G}_n)$  and
        \[
         \tau_n(s_1^{\prime},
        \dots, s_n^{\prime}) = (s_{1\tau}^{\prime}, \dots, s_{n\tau}^{\prime})
        = \pi_n\left(\beta\dfrac{z_1 - a}{1 - \bar{a}z_1}, \dots, \beta\dfrac{z_n - a}{1 - \bar{a}z_n}\right).
        \]
        Clearly for any $n \in \mathbb{N}$,
        \[
        (1-\bar{a}z_1)\dots(1-\bar{a}z_n)
        = 1 + \sum\limits_{i=1}^{n}(-\bar{a})^is_i
        \]
        and
        \[
        (z_1 - a)\dots(z_n - a)
        = \sum\limits_{i=0}^{n}(-a)^is_{n-i}\;,
        \]
        where $s_0 = 1$.
        Suppose the assertion is true for $n$. To complete the
        proof, we need to show that
\begin{align}\label{eqn:001}
& s_{i\tau} \notag \\ & =
\dfrac{\beta^i}{(1-\bar{a}z_1)\dots(1-\bar{a}z_{n+1})}\Bigg[\sum\limits_{j=1}^{i-1}(-a)^{i-j}
\bigg\{\sum\limits_{k=0}^{j}{{n+1-j} \choose {i-j+k}}{j \choose k}|a|^{2k}\bigg\}s_j \notag \\
             & \quad + \sum\limits_{j=i}^{n+1}(-\bar{a})^{j-i}\bigg\{\sum\limits_{k=0}^{i}
             {j \choose {i-k}}{{n+1-j} \choose k}|a|^{2k}\bigg\}s_j + (-a)^i{{n+1} \choose
             j}\Bigg].
\end{align}
This requires merely a few steps of routine calculations and we
skip it.

    \end{proof}
\end{lem}

\noindent The next result follows naturally and by Lemma \ref{lem1}, it provides transformation of a $\Gamma_n$-contraction under automorphisms of the symmetrized polydisc.

\begin{lem}\label{lem2}
    For $(S_1, \dots ,S_{n-1},P)$ and $\tau$ as above, $(S_{1\tau}, \dots ,S_{(n-1)\tau},P_{\tau})$ is a $\Gamma_n$-contraction.
    \begin{proof}
        We show that $\Gamma_n$ is a spectral set of
        $(S_{1\tau}, \dots ,S_{(n-1)\tau},P_{\tau})$. Let $f$ be a polynomial over $\mathbb{C}$ in $n$-variables. Then
        \begingroup
        \allowdisplaybreaks
        \begin{align*}
        \|f(S_{1\tau}, \dots ,S_{(n-1)\tau},P_{\tau})\| &= \|f\circ\tau(S_1, \dots ,S_{n-1},P)\|\\
        & \leq \|f\circ \tau\|_{\infty,\Gamma_n} = \sup\limits_{z\in \Gamma_n}|f(\tau(z))|
        \leq \|f\|_{\infty,\Gamma_n},
        \end{align*}
        \endgroup
         since $\tau(z)\in \Gamma_n$ for all $z \in \Gamma_n$ and
         hence  $(S_{1\tau}, \dots ,S_{(n-1)\tau},P_{\tau})$ is a $\Gamma_n$-contraction.
        \end{proof}

\end{lem}

\section{Unitary equivalence of two $\Gamma_n$-contractions}

\noindent In this section, we find a necessary and sufficient condition under which two $\Gamma_n$-contractions become unitarily equivalent. We begin with a result from the literature.

\begin{lem}[\cite{S:B}, Lemma 2.17]\label{lem3}
	For any commuting tuple $(S_1, \dots ,S_{n-1},P)$ of operators,
	$\sigma(S_1, \dots, S_{n-1}, P) \subseteq \Gamma_n$ if and only if
	\[
	I-\bar{a}S_1+\bar{a}^2S_2+ \dots +(-\bar{a})^{n-1}S_{n-1}+(-\bar{a})^nP
	\]
	is invertible for all $a \in \mathbb{D}$.
\end{lem}

 Suppose $(S_1, \dots, S_{n-1},P)$ is a $\Gamma_n$-contraction on a Hilbert space $\mathcal{H}$. For $a \in \mathbb C$, define $T_a = I + \sum\limits_{i=1}^{n-1}(-\bar{a})^iS_i + (-\bar{a})^nP$. Consider the set $\Lambda_{\Sigma} = \left\{ a \in \mathbb{C}: T_a \text{ is invertible}\right\}$.
By the above lemma it is clear that $\Lambda_{\Sigma} \supseteq \mathbb{D}$. Let us define 
\[
\Theta_{\Sigma}(a) = \left(S_{1\tau a}, \dots, S_{(n-1)\tau a}, S_{n\tau a}\right) \text{ for all } a \in \Lambda_{\Sigma},
\]
where
\begingroup
\allowdisplaybreaks
\begin{align*}
& S_{i\tau a}\\
= &   \left(\sum\limits_{j=1}^{i-1}(-a)^{i-j}
	\left[\sum\limits_{k=0}^j\left\{{{n-j}\choose{i-j+k}}{{j}\choose{k}}|a|^{2k}\right\}S_j\right]\right)T_a^{-1} \\
	 +& \left(\sum\limits_{j = i}^n(-\bar{a})^{j-i}\left[\sum\limits_{k=0}^i{{j}\choose {i-k}}{{n-j}\choose {k}}|a|^{2k}\right]S_j \right)T_a^{-1} + (-a)^i{{n}\choose{i}}T_a^{-1}.
\end{align*}
\endgroup
\begin{defn}
	Let $\Sigma = (S_1, \dots, S_{n-1},P)$ and $\Sigma^{\prime} = (S_1', \dots, S_{n-1}', P')$ be two $\Gamma_n$-contractions on $\mathcal{H}$ and $\mathcal{H}'$ respectively. Then we say that $\Theta_{\Sigma}$ and $\Theta_{\Sigma^{\prime}}$ coincide if for each $i = 1, \dots, n$ and $a\in \Lambda_{\Sigma}\cap \Lambda_{\Sigma^{\prime}}$, $S_{i\tau a}$ is unitarily equivalent to $S_{i\tau a}'$ by the same unitary.
\end{defn}
\begin{thm}\label{thm:equiv}
	Let $\Sigma = (S_1, \dots, S_{n-1},P)$ and $\Sigma^{\prime} = (S_1', \dots, S_{n-1}', P')$ be two $\Gamma_n$-contractions on $\mathcal{H}$ and $\mathcal{H}'$ respectively. Then $(S_1, \dots, S_{n-1},P)$ is unitarily equivalent to $(S_1', \dots, S_{n-1}', P')$ if and only if $\Theta_{\Sigma}$ and $\Theta_{\Sigma^{\prime}}$ coincide.
\end{thm}
\begin{proof}
	Suppose $\Theta_{\Sigma}$ and $\Theta_{\Sigma^{\prime}}$ coincide. Then by definition there exists a unitary $U : \mathcal{H} \rightarrow \mathcal{H}'$ such that $US_{i\tau a} = S_{i\tau a}'U$ for all $a\in \Lambda_{\Sigma}\cap \Lambda_{\Sigma^{\prime}}$ and $i= 1,\dots, n$. In particular, if $a=0$, then $\Theta_{\Sigma}(0) = \Sigma$ and $\Theta_{\Sigma^{\prime}}(0) = \Sigma^{\prime}$. Therefore, $US_i = S_i'U$ for all $i = 1,\dots, n-1$ and $UP = P'U$. Hence $\Sigma$ is unitarily equivalent to $\Sigma^{\prime}$. \\
	Conversely, suppose $\Sigma$ is unitarily equivalent to $\Sigma^{\prime}$. Then there exists a unitary $U : \mathcal{H} \rightarrow \mathcal{H}'$ such that $US_i = S_i'U$ and $UP = P'U$. To prove $US_{i\tau a}=S_{i\tau a}'U$ for all $a\in \Lambda_{\Sigma}\cap \Lambda_{\Sigma^{\prime}}$ it suffices to prove $UT_a^{-1} = T_a'^{-1}U$ for all $a\in \Lambda_{\Sigma}\cap \Lambda_{\Sigma^{\prime}}$. Now it is clear that $UT_a = T_a'U$ for all $a\in \Lambda_{\Sigma}\cap \Lambda_{\Sigma^{\prime}}$ and hence $UT_a^{-1} = T_a'^{-1}U$.
\end{proof}

\section{The automorphisms and the fundamental operator pair}

\noindent Before presenting the main result of this paper, we
recall few results from the literature which we shall use in the
proof the main theorem.

\begin{thm}[\cite{S:P}, Theorem 6.4] \label{thm3}
    If $(S_1,S_2,P)$ is a $\Gamma_3$-contraction,
    then $\left(\frac{1}{3}S_1+\frac{\beta}{3}S_2,\beta P \right)$
    is a $\Gamma$-contraction, for all $\beta \in \mathbb{T}$.
\end{thm}

\begin{thm}[{\cite{B:P}}, Theorem 2.2]\label{thm6}
    Suppose $(S_1,  \dots ,S_{n-1},P)$ is a $\Gamma_n$-contraction
    with commuting fundamental operator tuple $(A_1,  \dots ,A_{n-1})$. Then for each $i, j = 1,  \dots ,n-1$ we have
    \[
    S_i^*S_j - S_{n-j}^*S_{n-i} = D_P(A_i^*A_j - A_{n-j}^*A_{n-i})D_P.
    \]
\end{thm}

\noindent Now we state and prove the main result of this paper.

\begin{thm}\label{thm13}
    Let $(S_1,S_2,P)$ be a $\Gamma_3$-contraction with commuting fundamental operator pair $(A_1,A_2)$. Then
    
    \begin{enumerate}
    	\item[(i)] $\left(\bar{a}I-aA_2^*\right)\left(aI-\bar{a}A_2\right)\leq \left((1+|a|^2)I-aA_1^*\right)\left((1+|a|^2)I-\bar{a}A_1\right)$, for all $a\in \overline{\mathbb{D}}$\,;
    	\item[(ii)] if $r(A_1) \leq 2 $, then $\|(aI-\bar{a}A_2)\big((1+|a|^2)I-\bar{a}A_1\big)^{-1}\| \leq 1$,  for all $ a \in \mathbb{D}$.
    \end{enumerate}
    Moreover, if the inequality above is strict then for any $\tau =\tau_m \in Aut({\mathbb G}_3)$,
    there exists a unitary \;$\mathcal{U}:\mathcal{D}_{P_{\tau}} \rightarrow \mathcal{D}_P$
    such that the fundamental operator pair $(A_{1\tau},A_{2\tau})$ of $(S_{1\tau},S_{2\tau},P_{\tau})$ is given by
    \begingroup
    \allowdisplaybreaks
    \begin{align*}
    A_{1\tau} & = \mathcal{U}^*T^{-1/2}\beta \bigg( (1+3|a|^2)A_1 + a^2(3+|a|^2)A_1^* -
    2\bar{a}A_2 \\
    & \quad - 2a^3A_2^* - a(1+|a|^2)(A_1^*A_1 - A_2^*A_2) - 3a(1 + |a|^2)I\bigg)T^{-1/2}\mathcal{U}\;, \\
    A_{2\tau} & = \mathcal{U}^*T^{-1/2}\beta^2\bigg( -2(aA_1 + a^3 A_1^*) + A_2 + a^4 A_2^* \\ & + a^2(A_1^*A_1 - A_2^*A_2)
     + 3a^2I\bigg) T^{-1/2}\mathcal{U}\,,
     \end{align*}

where
\[
     T = (1+|a|^2 +|a|^4)I-(1+|a|^2)(aA_1^*+\bar{a}A_1) +
     (a^2A_2^*+\bar{a}^2A_2)\\
      + |a|^2(A_1^*A_1-A_2^*A_2)
\]
and $\tau=\tau_m$ as in $($\ref{eqn:01}$)\;\; \& \;\;
($\ref{eqn:02}$)$.
     \endgroup

\end{thm}

\begin{proof}

    We apply Lemma (\ref{lem1}) for $n=3$ and get
    \begingroup
    \allowdisplaybreaks
    \begin{align*}
    & \tau(s_1,s_2,p)\\
    = & \tau(z_1+z_2+z_3,z_1z_2+z_1z_3+z_2z_3,z_1z_2z_3)\\
    = & \left(\beta\sum\limits_{i=1}^{3}\dfrac{z_i - a}{1 - \bar{a}z_i},  \beta^2\sum\limits_{1 \leq i<j\leq 3}\dfrac{(z_i - a)(z_j - a)}{(1 - \bar{a}z_i)(1 - \bar{a}z_j)}, \beta^3\prod\limits_{i=1}^{3}\dfrac{(z_i - a)}{(1 - \bar{a}z_i)}\right)\\
    = & \bigg(\beta \dfrac{((1 + 2|a|^2)s_1 - \bar{a}(2 + |a|^2)s_2 + 3\bar{a}^2p - 3a)}{1 - \bar{a}s_1 + \bar{a}^2s_2 - \bar{a}^3p},\\
    &\; \beta^2 \dfrac{(-a(2 + |a|^2)s_1 + (1+2|a|^2)s_2-3\bar{a}p+3a^2)}{1-\bar{a}s_1+\bar{a}^2s_2-\bar{a}^3p},\beta^3\dfrac{(p-as_2+a^2s_1-a^3)}{1-\bar{a}s_1+\bar{a}^2s_2-\bar{a}^3p}\bigg).
    \end{align*}
    \endgroup
    It is obvious that $\tau$ can be defined on the open set
    \[
    \Gamma_{3a} = \left\{\left(\sum_{i=1}^{3}z_i,\sum_{1 \leq i<j\leq 3}z_iz_j,\prod_{i=1}^{3}z_i\right) : |z_i| < 1/|a|,  i=1,2,3\right\},
    \]

    \noindent which contains $\Gamma_3$. Clearly
    \[
    \left(S_{1\tau},S_{2\tau},P_{\tau}\right) = \tau(S_1,S_2,P),
    \]
     where
     \begingroup
     \allowdisplaybreaks
    \begin{align*}
    & S_{1\tau}  = \beta\big((1 + 2|a|^2)S_1 - \bar{a}(2 + |a|^2)S_2 + 3\bar{a}^2P - 3aI\big)\big(I - \bar{a}S_1 + \bar{a}^2S_2 - \bar{a}^3P\big)^{-1},\\
    & S_{2\tau} = \beta^2 \big(-a(2 + |a|^2)S_1 + (1 + 2|a|^2)S_2 - 3\bar{a}P + 3a^2I\big)\big(I - \bar{a}S_1 + \bar{a}^2S_2 - \bar{a}^3P\big)^{-1},\\
    & \qquad\qquad\qquad\qquad\qquad\qquad \text{and } \\
    & P_{\tau}  = \beta^3\big(P - aS_2 + a^2S_1 - a^3I\big)\big(I - \bar{a}S_1 + \bar{a}^2S_2 - \bar{a}^3P\big)^{-1}.
    \end{align*}
    \endgroup
     Here
    \begingroup
    \allowdisplaybreaks
    \begin{align*}
     & D_{P_{\tau}}^2 \\
     & = (I-P_{\tau}^*P_{\tau})\\
     & =I - \big(I - aS_1^* + a^2S_2^* - a^3P^*\big)^{-1}\big(P^* - \bar{a}S_2^* + \bar{a}^2S_1^* - \bar{a}^3I\big)\big(P - aS_2 \\
     & \quad + a^2S_1 - a^3I\big)\big(I - \bar{a}S_1 + \bar{a}^2S_2 - \bar{a}^3P\big)^{-1}\\
     & = \big(I - aS_1^* + a^2S_2^* - a^3P^*\big)^{-1} \Big\{\big(I - aS_1^* + a^2S_2^* - a^3P^*\big)\big(I - \bar{a}S_1 + \bar{a}^2S_2\\
     & \quad - \bar{a}^3P\big) - \big(P^* - \bar{a}S_2^* + \bar{a}^2S_1^* - \bar{a}^3I\big)\big(P - aS_2 + a^2S_1 - a^3I\big) \Big\}\big(I - \bar{a}S_1\\
     & \quad + \bar{a}^2S_2 - \bar{a}^3P\big)^{-1}\\
     &= (1 - |a|^2)\big(I - aS_1^* + a^2S_2^* - a^3P^*\big)^{-1}D_P \Big\{(1 + |a|^2 + |a|^4)I - (1+|a|^2)\\
     & \quad (aA_1^* + \bar{a}A_1) + (a^2A_2^* + \bar{a}^2A_2) + |a|^2(A_1^*A_1 - A_2^*A_2) \Big\}D_P\big(I - \bar{a}S_1 + \bar{a}^2S_2\\
     & \quad - \bar{a}^3P\big)^{-1}.
    \end{align*}
    \endgroup
     Now we show that the operator
    \[(1 + |a|^2 + |a|^4)I - (1+|a|^2)(aA_1^* + \bar{a}A_1) + (a^2A_2^* + \bar{a}^2A_2) + |a|^2(A_1^*A_1 - A_2^*A_2)\]
    defined on $\mathcal{D}_p$ is positive.
    Suppose $T = (1 + |a|^2 + |a|^4)I - (1+|a|^2)(aA_1^* + \bar{a}A_1) + (a^2A_2^* + \bar{a}^2A_2) + |a|^2(A_1^*A_1 - A_2^*A_2) $. Now for $h\in \mathcal{H}$ we have
    \begingroup
    \allowdisplaybreaks
    \begin{align*}
    &\langle TD_Ph,D_Ph\rangle\\
    & = \langle D_PTD_Ph,h\rangle\\
    & = \dfrac{1}{1 - |a|^2}\langle (I - aS_1^* + a^2S_2^* - a^3P^*)D_{P_{\tau}}^2(I - \bar{a}S_1 + \bar{a}^2S_2 - \bar{a}^3P)h,h \rangle\\
    & = \dfrac{1}{1 - |a|^2}\langle D_{P_{\tau}}^2(I - \bar{a}S_1 + \bar{a}^2S_2 - \bar{a}^3P)h,(I - \bar{a}S_1 + \bar{a}^2S_2 - \bar{a}^3P)h \rangle\\
    & = \dfrac{1}{1 - |a|^2}\langle D_{P_{\tau}}(I - \bar{a}S_1 + \bar{a}^2S_2 - \bar{a}^3P)h, D_{P_{\tau}}(I - \bar{a}S_1 + \bar{a}^2S_2 - \bar{a}^3P)h \rangle\\
    & = \dfrac{1}{1 - |a|^2}\|D_{P_{\tau}}(I - \bar{a}S_1 + \bar{a}^2S_2 - \bar{a}^3P)h\|^2\\
    &\geq 0.
    \end{align*}
    \endgroup
    Also
    \begingroup
    \allowdisplaybreaks
    \begin{align*}
    & (1 + |a|^2 + |a|^4)I - (1 + |a|^2)(aA_1^* + \bar{a} A_1) + (a^2 A_2^* + \bar{a}^2 A_2) + |a|^2(A_1^*A_1\\
    & - A_2^*A_2) \\
    & = \big((1 + |a|^2 )I - aA_1^* \big)\big((1 + |a|^2)I - \bar{a}A_1 \big) - \big(\bar{a}I - aA_2^* \big)\big(aI - \bar{a}A_2\big) \\
    & \geq 0.
    \end{align*}
    \endgroup
    This implies
    $$ \big(\bar{a}I - aA_2^* \big)\big(aI - \bar{a}A_2 \big) \leq \big((1 + |a|^2)I - aA_1^* \big)\big((1 + |a|^2)I - \bar{a}A_1 \big)$$ for all $ a \in \mathbb{D}$. By continuity we have that 
    \[
    \big(\bar{a}I - aA_2^* \big)\big(aI - \bar{a}A_2 \big) \leq \big((1 + |a|^2)I - aA_1^* \big)\big((1 + |a|^2)I - \bar{a}A_1 \big)
    \]
    for all $a \in \overline{\mathbb{D}}$.

   \noindent  \textbf{Claim:} If $r(A_1)\leq 2$, then $ \left((1 + |a|^2)I - \bar{a}A_1 \right)$ is invertible for all $a \in \mathbb{D}$.\\
    If $a = 0$ then it is obvious. If $0<|a| < 1$, then $\dfrac{1}{|a|} + |a| > 2.$
    Now
    \begingroup
    \allowdisplaybreaks
    \begin{align*}
         r\left(\dfrac{\bar{a}A_1}{1+|a|^2}\right)& =
         \lim\limits_{n\rightarrow \infty}\left\|\left(\dfrac{\bar{a}A_1}{1+|a|^2}\right)^n\right\|^{\frac{1}{n}}\\
        & =  \lim\limits_{n \rightarrow \infty}
        \left|\left(\dfrac{\bar{a}}{1+|a|^2}\right)^n\right|^\frac{1}{n}\|A_1^n\|^\frac{1}{n}\\
        & =\left|\dfrac{\bar{a}}{1+|a|^2}\right|\lim\limits_{n \rightarrow \infty}\|A_1^n\|^{\frac{1}{n}}\\
        & < 1.
    \end{align*}
    \endgroup
    Therefore, $ \left((1 + |a|^2)I - \bar{a}A_1\right) $ is invertible.
    This completes the proof of the claim. So we have
    \[
     \big((1 + |a|^2)I - aA_1^* \big)^{-1}\big(\bar{a}I-aA_2^* \big)
    \big(aI-\bar{a}A_2 \big)\big((1 + |a|^2)I - \bar{a}A_1 \big)^{-1} \leq I
    \] 
    and consequently
    \[
    \|\big(aI-\bar{a}A_2 \big)\big((1 + |a|^2)I - \bar{a}A_1 \big)^{-1}\| \leq
    1, \text{ for all } a \in \mathbb{D}\,.
    \]
Suppose  $\|\big(aI-\bar{a}A_2 \big)\big((1 + |a|^2)I - \bar{a}A_1
\big)^{-1}\| < 1$.
    Then $T$ is invertible. Let $X = (1-|a|^2)^{1/2}T^{1/2}D_P(I - \bar{a}S_1 + \bar{a}^2S_2 - \bar{a}^3P)^{-1}$.
    Then $X$ is an operator from $\mathcal{H}$ to $\mathcal{D}_P$. Also $D_{P_{\tau}}^2 = X^*X$
    and $\overline{\text{Ran}}X = \mathcal{D}_P$ as T is invertible. Now define
    $$ \mathcal{U} : \mathcal{D}_{P_{\tau}} \rightarrow \overline{\text{Ran}}X
    = \mathcal{D}_P $$ $$ D_{P_{\tau}}h \mapsto Xh. $$ Clearly $\mathcal{U}$ is onto.
    Moreover, $$ \|\mathcal{U} D_{P_{\tau}}h\|^2 = \|Xh\|^2 = \langle X^*Xh,h \rangle
    = \langle D_{P_{\tau}}^2h,h\rangle = \|D_{P_{\tau}}h\|^2. $$
    So $\mathcal{U}$ is a surjective isomery i.e., a unitary.
    Also
    \begingroup
    \allowdisplaybreaks
    \begin{align*}
    & S_{1\tau} - S_{2\tau}^*P_\tau\\
    & = \beta  \Big\{\big((1 + 2|a|^2)S_1 - \bar{a}(2 + |a|^2)S_2 + 3\bar{a}^2P - 3aI\big)\big(I - \bar{a}S_1 + \bar{a}^2S_2 -\\
    &\quad \bar{a}^3P\big)^{-1} - \big (I - aS_1^* + a^2S_2^* - a^3P^* \big )^{-1}\big (-\bar{a}(2+|a|^2)S_1^* + (1+2|a|^2)S_2^*\\
    &\quad - 3aP^* + 3\bar{a}^2I\big )\big (P - aS_2 + a^2S_1 - a^3I\big )\big (I - \bar{a}S_1 + \bar{a}^2S_2 - \bar{a}^3P\big )^{-1} \big ]    \\
    & = \beta \big (I - aS_1^* + a^2S_2^* - a^3P^* \big )^{-1}D_P  \Big\{(1+2|a|^2-3|a|^4)A_1 + a^2(3-2|a|^2\\
    &\quad -|a|^4)A_1^* - \bar{a}(2-2|a|^2)A_2 - a^3(2-2|a|^2)A_2^* - a(1-|a|^4)(A_1^*A_1-\\
    &\quad A_2^*A_2)-3a(1-|a|^4)I \Big\}D_P \big(I - \bar{a}S_1 + \bar{a}^2S_2 - \bar{a}^3P \big)^{-1}\\
    & = (1-|a|^2)\big (I - aS_1^* + a^2S_2^* - a^3P^* \big )^{-1} \beta D_P  \Big\{ (1+3|a|^2)A_1 + a^2(3+\\
    & \quad |a|^2)A_1^*- 2\bar{a}A_2 - 2a^3A_2^* - a(1+|a|^2)(A_1^*A_1-A_2^*A_2) - 3a(1+|a|^2)I  \Big\}\\
    & \quad  D_P\big(I - \bar{a}S_1+\bar{a}^2S_2 - \bar{a}^3P \big)^{-1}\\
    & = X^*T^{-1/2}\beta  \Big\{ (1+3|a|^2)A_1 + a^2(3+|a|^2)A_1^* - 2\bar{a}A_2 - 2a^3A_2^* - a(1+\\
    & \quad|a|^2)(A_1^*A_1-A_2^*A_2) - 3a(1+|a|^2)I  \Big\}T^{-1/2}X\\
    & = D_{P_{\tau}}\mathcal{U}^*T^{-1/2}\beta  \Big\{ (1+3|a|^2)A_1 + a^2(3+|a|^2)A_1^* - 2\bar{a}A_2 - 2a^3A_2^* - a(1\\
    & \quad+|a|^2)(A_1^*A_1-A_2^*A_2) - 3a(1+|a|^2)I  \Big\}T^{-1/2}\mathcal{U}D_{P_{\tau}}.
    \end{align*}
    \endgroup
    Again since $ S_{1\tau} - S_{2\tau}^*P_{\tau} = D_{P_{\tau}}A_{1\tau}D_{P_{\tau}}$ and $A_{1\tau}$ is unique, we have \\

    \noindent $ A_{1\tau} = \mathcal{U}^*T^{-1/2}\beta  \Big\{ (1+3|a|^2)A_1 + a^2(3+|a|^2)A_1^* - 2\bar{a}A_2 - 2a^3A_2^* - a(1+|a|^2)(A_1^*A_1-A_2^*A_2) - 3a(1+|a|^2)I  \Big\}T^{-1/2}\mathcal{U}. $ \\

    \noindent Again
    \begingroup
    \allowdisplaybreaks
    \begin{align*}
    &S_{2\tau} - S_{1\tau}^*P_{\tau}\\
    & =\beta^2 \Big\{ \big(-a(2 + |a|^2)S_1 + (1 + 2|a|^2)S_2 - 3\bar{a}P + 3a^2I\big)\big(I - \bar{a}S_1 + \bar{a}^2S_2\\
    &\quad -\bar{a}^3P\big)^{-1} - \big (I - aS_1^* + a^2S_2^* - a^3P^* \big )^{-1}\big( (1+|a|^2)S_1^*-a(2+|a|^2)S_2^*\\
    &\quad +3a^2P^*-3\bar{a}I \big)\big (P - aS_2 + a^2S_1 - a^3I\big )\big (I - \bar{a}S_1 + \bar{a}^2S_2 - \bar{a}^3P\big )^{-1}  \Big\} \\
    & =\beta^2 (1-|a|^2)\big (I - aS_1^* + a^2S_2^* - a^3P^* \big )^{-1}D_P  \Big\{ -2(aA_1+a^3A_1^*) + A_2 +\\
    & \quad a^4A_2^* + a^2(A_1^*A_1-A_2^*A_2) + 3a^2I \Big\}D_P\big (I - \bar{a}S_1 + \bar{a}^2S_2 - \bar{a}^3P\big )^{-1}\\
    & =D_{P_{\tau}}\mathcal{U}^*T^{-1/2}\beta^2  \Big\{ -2(aA_1+a^3A_1^*) + A_2 +a^4A_2^* + a^2(A_1^*A_1-A_2^*A_2) +\\
    &\quad 3a^2I \Big\}T^{-1/2}\mathcal{U}D_{P_{\tau}}.
    \end{align*}
    \endgroup
    \noindent Since $S_{2\tau} - S_{1\tau}^*P_{\tau} = D_{P_{\tau}}A_{2\tau}D_{P_{\tau}}$ and $A_{2\tau}$ is unique, we
    have that

    \begin{align*}
     A_{2\tau} = \mathcal{U}^*T^{-1/2}\beta^2 & \Big\{ -2(aA_1+a^3A_1^*) + A_2 +a^4A_2^* \\
     & \quad + a^2(A_1^*A_1-A_2^*A_2) + 3a^2I \Big\}T^{-1/2}\mathcal{U}.
     \end{align*}
     The proof is now complete.

\end{proof}

\section{Non-commuting fundamental operator tuple: an example}

\noindent In Theorem 3.6 in \cite{costara1}, Costara has shown that corresponding to every point $(s_1,\dots ,s_{n-1},p)\in \Gamma_n$, there is a unique point $(c_1,\dots ,c_{n-1})\in \Gamma_{n-1}$ such that
\[
s_i=c_i+\bar{c}_{n-i}p\,, \text{ for each } i=1,\dots , n-1.
\]
In the development of the theory of $\Gamma_n$-contraction, when we introduced the $\ft$-tuple associated with a $\Gamma_n$-contraction, it became clear that for a scalar $\Gamma_n$-contraction $(s_1,\dots ,s_{n-1},p)$ (which is nothing but a point in $\Gamma_n$), the $\ft$-tuple is just the $\Gamma_{n-1}$-contraction $(c_1,\dots, c_{n-1})$. Therefore, it is naturally asked whether the same result holds for a $\Gamma_n$-contraction.
Also, in the hypothesis of Theorem \ref{thm13}, we assumed the fact that the fundamental operator pair $(A_1,A_2)$ of the concerned $\Gamma_3$-contraction is commutative and that $r(A_1)\leq 2$. In this section, we shall show that the $\ft$-tuple of a $\Gamma_n$-contraction may not be a $\Gamma_{n-1}$-contraction, indeed it may not even be a commuting $n-1$ tuple. Here we provide a $\Gamma_3$-contraction $(S_1,S_2,P)$ on $\mathbb{C}^2$ that has non-commuting $\ft$-pair $(A_1,A_2)$ and that $r(A_1)>2$. Consider the following $2 \times 2$ matrix:
\[
A = 
\begin{pmatrix}
0 & \dfrac{3}{4}\\
1 & 0\\
\end{pmatrix}.
\]
Since $A$ is a contraction which dilates to a unitary $U$ and consequently we have that $\overline{\mathbb{D}}^3$ is a complete spectral set for $(A,A,A)$ which further implies that $\overline{\mathbb{D}}^3$ is a spectral set for $(A,A,A)$. Therefore, their symmetrization $(3A, 3A^2,A^3)$ is a $\Gamma_3$-contraction. Now consider
\[
(S_1, S_2, P) = \left( 
\begin{pmatrix}
0 & \dfrac{9}{4}\\
3 & 0\\
\end{pmatrix}, 
\begin{pmatrix}
\dfrac{9}{4} & 0\\
0 & \dfrac{9}{4}\\
\end{pmatrix},
\begin{pmatrix}
0 & \dfrac{9}{16} \\
\dfrac{3}{4} & 0 \\
\end{pmatrix}
\right).
\]
One can easily check that 
\begingroup
\allowdisplaybreaks
\begin{align*}
&D_P = (I-P^*P)^{1/2} =
\begin{pmatrix}
\dfrac{\sqrt{7}}{4} & 0\\
0 & \dfrac{\sqrt{175}}{16}
\end{pmatrix}\,,\\
&S_1 -S_2^*P  = \begin{pmatrix}
0 & \dfrac{63}{64}\\
\dfrac{21}{16} & 0
\end{pmatrix} = D_P \begin{pmatrix}
0 & \dfrac{9}{5}\\
\dfrac{12}{5} & 0
\end{pmatrix} D_P\\
& \qquad\qquad\qquad\qquad\quad\text{and }\\
&S_2 - S_1^*P  = \begin{pmatrix}
0 & 0\\
0 & \dfrac{63}{64}
\end{pmatrix} = D_P \begin{pmatrix}
0 & 0\\
0 & \dfrac{36}{25}
\end{pmatrix}D_P.
\end{align*}
\endgroup
This implies that \[
(A_1, A_2) = \left( \begin{pmatrix}
0 & \dfrac{9}{5}\\
\dfrac{12}{5} & 0
\end{pmatrix},\begin{pmatrix}
0 & 0\\
0 & \dfrac{36}{25}
\end{pmatrix}\right)
\]
is the $\ft$-pair of $(S_1, S_2,P)$. It is clear that $A_1A_2 \neq A_2A_1$ and $r(A_1) = \dfrac{\sqrt{108}}{5} >2$.


\end{document}